\documentclass[11pt]{article}
\usepackage{latexsym,amssymb}
\def\demo{\noindent{\bf Proof. }}
\def\QED{\hfill$\Box$}

\newtheorem{Theorem}{Theorem}[section]
\newtheorem{Lemma}[Theorem]{Lemma}
\newtheorem{Corollary}[Theorem]{Corollary}
\newtheorem{Proposition}[Theorem]{Proposition}

\newtheorem{Definition}[Theorem]{Definition}

\begin{document}

\topmargin3mm
\hoffset=-1cm
\voffset=-1.5cm
\begin{center}

{\Large\bf Rees cones and monomial rings of matroids}\\

\vspace{6mm}


\addtocounter{footnote}{1}


\medskip

Rafael
H. Villarreal\footnote{Partially supported
by CONACyT grant 49251-F and SNI, M\'exico.}\\

{\small Departamento de
Matem\'aticas}\vspace{-1mm}\\
{\small Centro de Investigaci\'on y de Estudios
Avanzados del
IPN}\vspace{-1mm}\\
{\small Apartado Postal
14--740}\vspace{-1mm}\\
{\small 07000 Mexico City, D.F.}\vspace{-1mm}\\
{\small
e-mail: {\tt
vila@math.cinvestav.mx}}\vspace{4mm}

\end{center}

\date{}

\begin{abstract}
\noindent Using linear algebra methods we study certain algebraic
properties 
of monomial rings and matroids. Let $I$ be a monomial ideal in 
a polynomial ring over an arbitrary field. If the Rees cone of $I$ 
is quasi-ideal, we express the normalization of the Rees algebra 
of $I$ in terms of an Ehrhart ring. We introduce the basis Rees cone 
of a matroid (or a polymatroid) and 
study their facets. Some applications to Rees algebras are presented.
It is shown that the basis monomial ideal of a matroid (or a
polymatroid) is normal. 
\end{abstract}

\section{Introduction}\label{Int}
Let $R=k[x_1,\ldots,x_n]$ be a polynomial ring over a field $k$ 
and let $I=(x^{v_1},\ldots,x^{v_q})$ be a proper monomial ideal 
of $R$ generated by $F=\{x^{v_1},\ldots,x^{v_q}\}$. As usual for 
$a=(a_i)$ in $\mathbb{N}^n$ we set $x^a=x_1^{a_1}\cdots x_n^{a_n}$. 
The {\it Rees cone\/} of $I$ is the rational polyhedral 
cone in $\mathbb{R}^{n+1}$, denoted by $\mathbb{R}_+{\cal A}'$,
consisting  
of all the non-negative linear combinations of the set 
$${\cal
A}':=\{e_1,\ldots,e_n,(v_1,1),\ldots,(v_q,1)\}\subset\mathbb{R}^{n+1},$$
where $e_i$ is the $i${\it th} unit vector. 
Rees cones have been used in commutative algebra to
study algebraic properties and invariants of blowup algebras of 
monomial ideals \cite{adrian-ainv,normali,wolmer65}. 
Notice that $\mathbb{R}_+{\cal A}'$ has dimension
$n+1$. Hence, by the finite basis theorem, the Rees cone is
defined by a finite system of linear inequalities, i.e., there is a 
unique (up to permutation of rows) integer matrix $B=(b_{ij})$ of
order $m\times(n+1)$ with
non-zero rows such that 
$$\mathbb{R}_+{\cal A}'=\{y\vert\, By\geq 0\},$$
the non-zero entries of each row of $B$ are relatively
prime, and none of the rows of $B$ can be omitted. In this situation
it is well known that the {\it facets\/} (faces of maximal dimension)
of $\mathbb{R}_+{\cal A}'$ are given by 
$$E_i=\{y\vert\,\textstyle\sum_{j}b_{ij}y_j=0\}\cap
\mathbb{R}_+{\cal A}'$$ 
for $i=1,\ldots,m$. Consider the convex polytope  $P={\rm
conv}(v_1,\ldots,v_q)$, where ``conv'' stands for ``convex hull''. The
{\it Ehrhart ring\/} of $P$ and the {\it Rees algebra\/} of $I$ 
are defined as the $k$-subrings of $k[x_1,\ldots,x_n,t]$ given by
$$
A(P)=k[\{x^at^b\vert\, a\in bP\cap\mathbb{Z}^n\}]\ \mbox{ and  }\ 
R[It]=k[x_1,\ldots,x_n,x^{v_1}t,\ldots,x^{v_q}t]
$$
respectively, where $t$ is a new variable. The {\it normalization\/} 
of $R[It]$, denoted by $\overline{R[It]}$, 
is the integral closure of $R[It]$ in its field of fractions. Notice
that the subgroup of
$\mathbb{Z}^{n+1}$ generated by ${\cal A}'$ is equal to
$\mathbb{Z}^{n+1}$. Thus using 
\cite[Theorem~7.2.28]{monalg}  we obtain: 
\begin{equation}\label{may6-06}
\overline{R[It]}=k[\{x^at^b\vert\, 
(a,b)\in\mathbb{R}_+{\cal A}'\cap\mathbb{Z}^{n+1}\}]\subset R[t].
\end{equation}
This formula explains the importance of Rees cones in the study 
of normalizations of Rees algebras. These algebras are important 
objects of study in algebra and geometry \cite{BSV,Vas}. Since 
we have the inclusion $A(P)[x_1,\ldots,x_n]\subset \overline{R[It]}$,
it is natural to ask 
when the equality occurs. We are able to prove equality if
$\mathbb{R}_+{\cal A}'$ is quasi-ideal (Theorem~\ref{nov27-05}), 
i.e., if the entries of $B$ are sufficiently nice or more precisely if
the first 
$n$ columns of $B$ have all their entries in $\{0,1\}$. 
If $I$ is square-free, in \cite{normali} it is shown 
that $\overline{R[It]}$ is equal to the symbolic Rees algebra of $I$
if and only if $\mathbb{R}_+{\cal A}'$ is {\it ideal\/}, i.e., if and
only if $B$ has all its entries in $\{0,\pm 1\}$. The terms ideal and
quasi-ideal refer to the facet structure of $\mathbb{R}_+{\cal A}'$
and not to the algebraic structure of an ideal in the sense of ring
theory. This terminology is in harmony with the combinatorial
optimization notion of {\it ideal clutter\/} and {\it ideal matrix\/}
introduced in \cite{cornu-book}. Indeed, let $A$ be the matrix
with column vectors $v_1,\ldots,v_q$ and assume that the entries of
$A$ are in $\{0,1\}$, by \cite[Theorem~3.2]{clutters} the Rees
cone $\mathbb{R}_+{\cal A}'$ is ideal if and only if 
the polyhedron $\{y\vert\, yA\geq{1}; y\geq{0}\}$ has
only integral vertices. 

Let $\mathbb{N}{\cal A}'$
be the 
subsemigroup of 
$\mathbb{N}^{n+1}$ generated by ${\cal A}'$, consisting of the 
linear combinations of ${\cal A}'$ with non-negative integer
coefficients. 
The Rees algebra of $I$ can be rewritten as
\begin{equation}\label{may6-06-2}
R[It]=
k[\{x^at^b\vert\, (a,b)\in\mathbb{N}{\cal A}'\}].
\end{equation}
Hence, by Eqs.~(\ref{may6-06}) and (\ref{may6-06-2}), we get that
$R[It]$ is a normal domain if and 
only if  the following equality holds: 
$$
\mathbb{N}{\cal A}'=
\mathbb{Z}^{n+1}\cap\mathbb{R}_+{\cal A}'.
$$
In geometric terms this equality means that $R[It]=\overline{R[It]}$ if
and only if ${\cal A}'$ is an integral Hilbert basis. As an
application we prove that $R[It]$ is normal if the Rees cone of $I$
is quasi-ideal 
and $k[Ft]=A(P)$ (Corollary~\ref{nov28-05}). There is a partial
converse of this result \cite{ehrhart,ITG}, namely, if all monomials
of $F$ have the same 
degree and $R[It]$ is normal, then $k[Ft]=A(P)$. 
  
\paragraph{\it Ideals and Rees cones of Matroids}  

Let $[n]=\{1,\ldots,n\}$ and let $M=([n],{\cal I})$ be a {\it matroid\/}
 on $[n]$, i.e., there is a collection $\cal I$ of 
subsets of $[n]$ containing the empty set and 
satisfying the following two conditions:
\begin{description}
\item{$(\mathrm{i}_1)$} If $I\in\cal I$ and $I'\subset I$, 
then $I'\in
\cal I$.
\item{$(\mathrm{i}_2)$} If $I_1$ and $I_2$ are in $\cal I$ and 
$|I_1|<|I_2|$, then there is an element $e$ of $I_2\setminus I_1$ such
that $I_1\cup\{e\}\in \cal I$.
\end{description} 

The members of $\cal I$ are the {\it independent sets\/} of $M$. A
maximal independent 
set of $M$ with respect to inclusion is called a {\it
basis}. The reader is referred to \cite{oxley} for the 
general theory of matroids. It is well known that all bases of $M$ have
the same number of elements, this common number $d$ is called the {\it
rank\/} of $M$. Let ${\cal B}=\{B_1,\ldots,B_q\}$ be the collection of
bases of $M$. For use below we set
$$
{\cal A}'=\{e_1,\ldots,e_n,f_{B_1},\ldots,f_{B_q}\}; \ \ \ \ 
f_{B_i}=e_{n+1}+\sum_{j\in B_i}e_j.
$$
 
An important aim here is to study the structure of the facets of
$\mathbb{R}_+{\cal A}'$,  
the {\it basis Rees cone\/} of $M$. An earlier result
 shows that the minimal primes of $I$ are in one to one correspondence 
with a certain collection of facets of $\mathbb{R}_+{\cal A}'$, see
\cite{normali}. 
One of our main results shows that the basis Rees cone of $M$ is
quasi-ideal (Theorem~\ref{sept21-05}). 
Let $F_M$ be the set of 
all monomials $x_{i_1}\cdots x_{i_d}$ in $R$ such 
that $\{{i_1},\ldots,{i_d}\}$ is a basis of $M$, and let $I=(F_M)$ 
be the {\it basis monomial ideal\/} of $M$.  
By a classical result of 
White \cite{white}, the subring $k[F_M]\subset R$ is normal. As an
application to Rees 
algebras, using  White's  result and the quasi-ideal property of the
basis 
Rees cone, we prove that $R[It]$ is normal
(Corollary~\ref{nwhite-gene}). At the 
end we introduce the family of bases of a discrete polymatroid and
show how to  
generalize Theorem~\ref{sept21-05} and Corollary~\ref{nwhite-gene} to
discrete polymatroids.

Along the paper we introduce most of the 
notions that are relevant for our purposes. 
For unexplained
terminology and notation on polyhedral geometry we refer to
\cite{Bron,webster}.

\section{Quasi-ideal Rees cones} 

Let $R=k[x_1,\ldots,x_n]$ be a polynomial ring over a field $k$ 
and let $I=(x^{v_1},\ldots,x^{v_q})$ be a proper monomial ideal 
of $R$ generated by $F=\{x^{v_1},\ldots,x^{v_q}\}$. 
The {\it Rees cone\/} of $I$ is the polyhedral 
cone in $\mathbb{R}^{n+1}$, denoted by $\mathbb{R}_+{\cal A}'$,
generated by the set 
$${\cal
A}':=\{e_1,\ldots,e_n,(v_1,1),\ldots,(v_q,1)\}\subset\mathbb{R}^{n+1},$$  
where $e_i$ is the $i${\it th} unit vector.  See
\cite[Section~3]{normali} for information 
about some of the interesting properties of Rees cones of square-free 
monomial ideals. Consider 
the index set
$${\cal J}=\{i\,| \; 1\leq i\leq n \mbox{ and }\langle e_i,v_j\rangle=0
\mbox{ for 
some } j\}\cup \{n+1\},
$$
where $\langle\ ,\, \rangle$ denotes the standard 
inner product. Notice that the cone $\mathbb{R}_+{\cal A}'$ has dimension
$n+1$, i.e., the linear space generated by ${\cal A}'$ is equal to
$\mathbb{R}^{n+1}$. Hence using
\cite[Theorem~3.2.1]{webster} it is seen that 
 the Rees cone has a unique irreducible representation:
\begin{equation}\label{okayama-car} 
{\mathbb R}_+{\cal A}'=\left(\bigcap_{i\in \cal
J}H_{e_i}^+\right)\bigcap\left(  
\bigcap_{i=1}^r H_{\ell_i}^+\right)
\end{equation}
such that none of the closed halfspaces can be omitted from
Eq.~(\ref{okayama-car}), $0\neq \ell_i\in\mathbb{Z}^{n+1}$ for all
$i$, 
and the non-zero entries of
each $\ell_i$ are relatively prime. Here $H_{a}^+$ denotes 
the {\it closed halfspace\/} $H_a^+=\{y\vert\, \langle
y,a\rangle\geq 0\}$ and $H_a$ is the hyperplane through the 
origin with normal vector $a$. It is easy to see that the first $n$ 
entries of each $\ell_i$ are non-negative and that the last entry of each
$\ell_i$ is negative, the second assertion follows from the
irreducibility of Eq.~(\ref{okayama-car}).   

To avoid repetitions in this section we shall use the notation 
introduced in Section~\ref{Int}. 
Thus $A(P)$ denotes the Ehrhart ring of the convex 
polytope  $P$, $R[It]$ denotes the Rees algebra of $I$, and
$\overline{R[It]}$ denotes the normalization of $R[It]$.

\medskip

\noindent{\it Notation\/} For use below we set $[n]=\{1,\ldots,n\}$. 
\begin{Theorem}\label{nov27-05} 
If each vector $\ell_k$ has the form 
$\ell_k=-d_ke_{n+1}+\sum_{i\in A_k}e_i$ for some $A_k\subset[n]$ and
some $d_k\in\mathbb{N}$, then $A(P)[x_1,\ldots,x_n]
=\overline{R[It]}$.
\end{Theorem}

\demo Clearly the left hand side is contained in the right hand side
because $A(P)\subset\overline{R[It]}$. To prove the reverse 
inclusion let $x^at^b=x_1^{a_1}\cdots
x_n^{a_n}t^b\in\overline{R[It]}$ be a minimal generator,  
that is, $(a,b)$ cannot be written as a sum of two non-zero integral
vectors in the Rees cone  
$\mathbb{R}_+{\cal A}'$. If $b=0$, then 
$x^at^b=x_i$ for some $i\in[n]$. Thus we may assume that $a_i\geq 1$ 
for $i\leq m$, $a_i=0$ for $i>m$, and $b\geq 1$. 

Case (I): $\langle(a,b),\ell_i\rangle>0$ for all $i$. The vector
$\gamma=(a,b)-e_1$  
satisfies $\langle\gamma,\ell_i\rangle\geq 0$ for all $i$, that is, 
$\gamma\in\mathbb{R}_+{\cal A}'$. Hence, since $(a,b)=e_1+\gamma$ and
$\gamma\neq 0$, we derive a contradiction. 

Case (II): $\langle(a,b),\ell_i\rangle=0$ for some $i$. We may assume
that
$$
\{\ell_j\vert\, \langle(a,b),\ell_j\rangle=0\}=\{\ell_1,\ldots,\ell_p\}.
$$

Subcase (II.a): $e_i\in H_{\ell_1}\cap\cdots\cap H_{\ell_p}$ for some
$i\in[m]=\{1,\ldots,m\}$.  
It is not hard to verify that the vector 
$\gamma=(a,b)-e_i$ satisfies $\langle\gamma,\ell_k\rangle\geq 0$ for
all $\ell_k$. Thus $\gamma\in\mathbb{R}_+{\cal A}'$, a contradiction because 
$(a,b)=e_i+\gamma$ and $\gamma\neq 0$. 

Subcase (II.b): $e_i\notin H_{\ell_1}\cap\cdots\cap H_{\ell_p}$ for
all $i\in[m]$. Since 
the vector $(a,b)$ belongs to the Rees cone it follows that we can write
$$
(a,b)=\lambda_1(v_1,1)+\cdots+\lambda_q(v_q,1) \ \ \ \ (\lambda_i\geq 0).
$$
Hence $a/b\in P$ and $a\in bP\cap\mathbb{Z}^n$, i.e., $x^at^b\in A(P)$. \QED

\begin{Definition}\rm If $\ell_1,\ldots,\ell_r$ satisfy the 
condition of Theorem~\ref{nov27-05} (resp. $d_k=1$ for all $k$), 
we say that the Rees cone 
of $I$ is {\it quasi-ideal} (resp. {\it ideal\/}).
\end{Definition}

As pointed out in the introduction the notions of ideal and
quasi-ideal Rees cones refer to the facet structure of
the polyhedral cone $\mathbb{R}_+{\cal A}'$, i.e, they refer to the 
irreducible representation of the Rees cone of $I$.

\begin{Corollary}\label{nov28-05} If the Rees cone of $I$ is quasi-ideal and 
$k[Ft]=A(P)$, then the Rees algebra $R[It]$ is normal.
\end{Corollary}

\demo By Theorem~\ref{nov27-05} we get 
$$
R[It]=k[Ft][x_1,\ldots,x_n]=A(P)[x_1,\ldots,x_n]=\overline{R[It]}.
\eqno\Box
$$

\section{Facets of Rees cones of matroids} 

It is worth mentioning that the following {\it exchange property\/}
allows to define 
the notion of a basis of a discrete polymatroid \cite{polymatroid}.
Later we will 
introduce this notion.

\begin{Theorem}{\rm\cite[Corollary~1.2.5, p.~18]{oxley}}
\label{exchangep} Let $\cal B$ be a non-empty family of 
subsets of $[n]$.
Then $\cal B$ is the collection of bases of a matroid
on $[n]$ if and
only if the following {exchange property} 
is satisfied: If $B_1$ and $B_2$ are members of $\cal B$ and $b_1\in
B_1\setminus B_2$, then there is an element $b_2\in B_2\setminus B_1$
such that $(B_1\setminus\{b_1\})\cup \{b_2\}$ is in $\cal
B$.
\end{Theorem}

The family of bases of a matroid satisfies the following {\it
symmetric exchange property\/} that will be used below (see \cite{kung}):

\begin{Theorem}\label{symmexchange} 
If $B_1$ and $B_2$ are bases of a matroid $M$ and $b_1\in
B_1\setminus B_2$, then there is an element $b_2\in B_2\setminus B_1$
such that both $(B_1\setminus\{b_1\})\cup \{b_2\}$ and
$(B_2\setminus\{b_2\})\cup \{b_1\}$ are bases of $M$.
\end{Theorem}

Let $M$ be a matroid of rank $d$ on the set $[n]=\{1,\ldots,n\}$ and
let ${\cal B}=\{B_1,\ldots,B_q\}$ be the collection of
bases of $M$. Notice that $|B_i|=d$ for all
$i$ and that $\cal B$ is a clutter on $[n]$, i.e., $B_i\not\subset B_j$
for $i\neq j$. 
For use below we set
$$
{\cal A}'=\{e_1,\ldots,e_n,f_{B_1},\ldots,f_{B_q}\}; \ \ \ \ 
f_{B_i}=e_{n+1}+\sum_{j\in B_i}e_j.
$$
We call $\mathbb{R}_+{\cal A}'$ the {\it basis Rees cone\/} of $M$.

\begin{Definition}\rm 
An integral matrix $C$ is 
called {\it totally unimodular\/} if each $i\times i$ 
minor of $C$ is $0$ or $\pm 1$ for 
all $i\geq 1$.
\end{Definition}

\begin{Lemma}\label{nov18-06} 
If $n=q\geq 2$ and $B_i=\{i\}$ for $i\in [n]$, 
then 
$$
\mathbb{R}_+{\cal A}'=H_{e_1}^+\cap\cdots\cap H_{e_{n+1}}^+\cap
H_{b}^+,
$$
where $b=(1,\ldots,1,-1)$. If $n=q=1$, then 
$\mathbb{R}_+{\cal A}'=H_{e_1}^+\cap H_{(1,-1)}^+$. Moreover these 
are irreducible representations of the Rees cone of $M$.
\end{Lemma}

\demo It follows readily using that the matrix whose columns are the vectors
in ${\cal A}'$ is totally unimodular. The lemma also follows at once 
from the expression for the irreducible representation of a Rees cone
given in  \cite[Theorem~3.2]{clutters}. \QED

\begin{Theorem}\label{sept21-05} 
If $E$ is a facet of $\mathbb{R}_+{\cal A}'$, then $E$ has the 
form $\mathbb{R}_+{\cal A}'\cap H_b$, 
where $b=e_i$ for some $i\in[n+1]$ or $b=(b_i)$ is an integral
vector 
such that $b_i\in\{0,1\}$ for 
$i\in[n]$ and $b_{n+1}\geq -d$. In particular the Rees cone of
$M$ 
is quasi-ideal.
\end{Theorem}

\demo The proof is by induction on $d$. The case $d=1$ is clear 
because of Lemma~\ref{nov18-06}. Assume $d\geq 2$. 
Let $E$ be a facet of $\mathbb{R}_+{\cal A}'$. Since the dimension of
the Rees cone is $n+1$, there is a unique 
integral vector $0\neq b=(b_i)$ whose non-zero entries are relatively
prime such that 
\begin{description}
\item{\rm (a)} $E=\mathbb{R}_+{\cal A}'\cap H_b\neq \mathbb{R}_+{\cal
A}'$,\ \,  $\mathbb{R}_+{\cal 
A}'\subset H_b^+$, and 
\item{\rm (b)} there is a linearly independent set 
${\cal D}\subset H_b\cap {\cal A}'$ with $|{\cal D}|=n$. 
\end{description}

If ${\cal D}=\{e_1,\ldots,e_n\}$ (resp. ${\cal D}\subset
\{f_{B_1},\ldots,f_{B_q}\}$), then $b=e_{n+1}$ 
(resp. $b$ is the vector $(1,\ldots,1,-d)$). The first assertion is
clear. To show the second assertion observe that the vectors $b$
and $(1,\ldots,1,-d)$ are in the orthogonal complement of 
$\mathbb{R}\{f_{B_1},\ldots,f_{B_q}\}$, and consequently they must be
equal because the non-zero entries of each vector are relatively
prime. Notice that $b_i\geq 0$ 
for $i=1,\ldots,n$. Thus we may assume that 
${\cal D}$ is the set $\{e_1,\ldots,e_s,f_{B_1},\ldots,f_{B_t}\}$, 
where  $s\in[n-1]$, $s+t=n$, and $\{1,\ldots,s\}$ is the set 
of all $i\in[n]$ such that $e_i\in H_b$. 

We may
assume that $1\in B_k$ for some $k\in [q]$. Indeed if 
$1$ is not in $\cup_{i=1}^q B_i$, then the facets of 
$\mathbb{R}_+{\cal A}'$ different from $H_{e_1}\cap \mathbb{R}_+{\cal
A}'$ 
are in one to one correspondence with the
facets of the Rees cone
$$
\mathbb{R}_+\{e_2,\ldots,e_n,f_{B_1},\ldots,f_{B_q}\}.
$$
Observe that this argument can be applied replacing $1$ by
any other element $j$ in
$[n]$ and $j$ not in $\cup_{i=1}^q B_i$. Hence we may as well assume
that the set $[n]$ is equal to $\cup_{i=1}^q B_i$.  

Assume that $1\notin B_1$.
Since $1\in B_k\setminus B_1$ for some $k$, 
by the symmetric exchange property of $\cal B$ 
(Theorem~\ref{symmexchange}), there is $j\in B_1\setminus B_k$ such 
that $(B_1\setminus\{j\})\cup\{1\}=B_i$ for some basis $B_i$. Hence 
$$
f_{B_1}+e_1=f_{B_i}+e_j\ \Rightarrow\ \langle f_{B_i},b\rangle=-
\langle e_j,b\rangle \ \Rightarrow\ \langle f_{B_i},b\rangle=
\langle e_j,b\rangle=0,
$$
i.e., $f_{B_i}$ and $e_j$ belong to $H_b$. Thus, as $1\in B_i$, by
replacing $f_{B_1}$ by $f_{B_i}$ in $\cal D$, 
we may assume from the outset that $1\in B_1$.  
Next assume that $t\geq 2$ and $1\notin B_2$. Since $1\in B_1\setminus B_2$, 
by the exchange property there is $j\in B_2\setminus B_1$ such 
that $(B_2\setminus\{j\})\cup\{1\}=B_i$ for some basis $B_i$. Hence 
$$
f_{B_2}+e_1=f_{B_i}+e_j\ \Rightarrow\ \langle f_{B_i},b\rangle=
\langle e_j,b\rangle=0,
$$
i.e., $f_{B_i}$ and $e_j$ belong to $H_b$. Notice that $i\neq 1$. Indeed
if $i=1$, then from the equality above we get that 
$\{f_{B_2},e_1,f_{B_1},e_j\}$ is linearly dependent, a contradiction
because this set is contained in $\cal D$. Thus, as $1\in B_i$, by
replacing $f_{B_2}$ by $f_{B_i}$ in $\cal D$, 
we may assume from the outset that $1\in B_1\cap B_2$.
Applying the arguments above repeatedly shows that we may assume that $1$ 
belongs to $B_i$ for $i=1,\ldots,t$. We may also assume that 
$B_1,\ldots,B_r$ is the set of all basis of $M$ such that $1\in B_i$,
where $t\leq r$. For $i\in [r]$, we set $C_i=B_i\setminus\{1\}$.
Notice that there is a matroid $M'$ on $[n]$ of rank $d-1$ whose collection of 
bases is $\{C_1,\ldots,C_r\}$. Consider the Rees cone
$\mathbb{R}_+{\cal A}''$ generated by 
$$
{\cal A}''=\{e_1,e_2,\ldots,e_n,f_{C_1},\ldots,f_{C_r}\}.
$$
Notice that $\mathbb{R}_+{\cal A}''\subset
H_b^+$ and that $\{e_1,e_2,\ldots,e_s,f_{C_1},\ldots,f_{C_t}\}$ is a
linearly independent set. Thus $\mathbb{R}_+{\cal A}''\cap H_b$ is a
facet of the cone $\mathbb{R}_+{\cal A}''$ and the result follows by
induction. \QED 

\medskip

The facets of the Rees cone of the ideal generated by all square-free
monomials of degree $d$ of $R$ were computed in 
\cite[Theorem~3.1]{adrian-ainv}, the result above can be seen as a 
generalization of this result.

\medskip

The set of all monomials $x_{i_1}\cdots x_{i_d}\in R$ such 
that $\{{i_1},\ldots,{i_d}\}$ is a basis of $M$ will be denoted by 
$F_M$ and the 
subsemigroup  generated by $F_M$ will be denoted 
by $\mathbb{M}_M$. The {\it basis monomial ring\/}
 of $M$ is the subring $k[F_M]\subset R$. The square-free monomial ideal 
$I({\cal B}):=(F_M)$ is called the {\it basis monomial ideal} of $M$. 

\begin{Proposition}[\rm\cite{white}]\label{nwhite} If 
$x^a$ is a monomial of degree $\ell{d}$ for some 
$\ell\in\mathbb{N}$ such that $(x^{a})^p\in\mathbb{M}_M$ for some 
$p\in\mathbb{N}\setminus\{0\}$, then $x^a\in\mathbb{M}_M$. In 
particular $k[F_M]$ is normal. 
\end{Proposition}

The next result is just a reinterpretation of
Proposition~\ref{nwhite}, which is adequate to examine 
the normality of the Rees algebra of $I({\cal B})$.  

\begin{Corollary}\label{nwhite-coro} Let $P$ be the convex hull 
of the set of all vectors $e_{i_1}+\cdots+e_{i_d}$ such that
$\{{i_1},\ldots,{i_d}\}$ is a 
basis of $M$. Then $A(P)=k[F_Mt]$. 
\end{Corollary}

\demo It suffices to show the inclusion $A(P)\subset k[F_Mt]$. 
Take $x^at^b\in A(P)$, that is, $a\in\mathbb{Z}^n\cap bP$. Hence 
$x^a$ has  
degree $bd$ and $(x^{a})^p\in\mathbb{M}_M$ for some positive integer
$p$. By  
Proposition~\ref{nwhite} we get $x^a\in\mathbb{M}_M$. It is rapidly
seen that  
$x^at^b$ is in $k[F_Mt]$, as required. \QED

\medskip

As an application to Rees algebras we obtain the following: 

\begin{Corollary}\label{nwhite-gene}
If $I=I({\cal B})$ is the basis monomial ideal of a
matroid, then $R[It]$ is normal. 
\end{Corollary}

\demo It follows from Corollaries~\ref{nov28-05} and
\ref{nwhite-coro} together with  Theorem~\ref{sept21-05}.  \QED 

\paragraph{Polymatroidal sets of monomials} For 
$a=(a_1,\ldots,a_n)\in\mathbb{N}^n$ we 
set $|a|=a_1+\cdots+a_n$. Let ${\cal A}=
\{v_1,\ldots,v_q\}\subset\mathbb{N}^n$ be the set of bases of a  
discrete polymatroid of rank $d$, i.e., $|v_i|=d$ for all $i$ 
and the following condition is satisfied:  given any two
$a=(a_i),c=(c_i)$ in ${\cal A}$, if $a_i> c_i$ for some index $i$,  
then there is an
index $j$ with $a_j < c_j$ such that $a-e_i+e_j$ is in $\cal A$.  
The set $F=\{{x}^{v_1},\ldots,{x}^{v_q}\}$ (resp. the ideal 
$I=(F)\subset R$) is called a 
{\it polymatroidal set of monomials\/} (resp. a 
{\it polymatroidal ideal\/}). Notice that the basis monomial ideal of 
a matroid is a polymatroidal ideal. We 
refer to \cite{polymatroid,fiber-type} for the theory of 
discrete polymatroids and polymatroidal ideals. Below we indicate how
to  
generalize Theorem~\ref{sept21-05} and Corollary~\ref{nwhite-gene} to 
discrete polymatroids.

\begin{Lemma}\label{dec15-05} The set 
$F'=\{x^{v_i}/x_1\colon\,  x_1\mbox{ occurs in }x^{v_i}\}$ is also
polymatroidal.  
\end{Lemma}

\begin{Lemma}\label{dec17-05} Let
$G=\{x^{u_1},\ldots,x^{u_s}\}$ and let $d=\max\{|u_i|:\, i\in[s]\}$.
Suppose that $F=\{x^{u_1},\ldots,x^{u_t}\}$ is 
the set of all $x^{u_i}$ of degree $d$. If 
 $Q={\rm conv}(u_1,\ldots,u_s)$ and $k[Gt]=A(Q)$, then $k[Ft]=A(P)$, 
where $P$ is the convex hull of ${u_1},\ldots,{u_t}$.
\end{Lemma}

\begin{Proposition}\label{nwhite-gene-poly} 
If $I=(F)$ is a polymatroidal ideal, then $R[It]$ is normal. 
\end{Proposition}

\demo Let ${\cal A}'=\{e_1,\ldots,e_n,(v_1,1),\ldots,(v_q,1)\}$ and 
let $\mathbb{R}_+{\cal A}'$ be the Rees cone of $I$. Using the 
proof of Theorem~\ref{sept21-05} together with 
Lemma~\ref{dec15-05} and the fact that $\cal A$ satisfies the 
symmetric exchange property \cite[Theorem~4.1, p.~241]{polymatroid}
it is not hard to see that the Rees cone $\mathbb{R}_+{\cal A}'$ is
quasi-ideal. By \cite[Theorem~6.1]{polymatroid} and
Lemma~\ref{dec17-05} we get $k[Ft]=A(P)$. Thus applying 
Corollary~\ref{nov28-05} we obtain that $R[It]$ is normal. \QED

\bibliographystyle{plain}

\end{document}